\documentclass{article}
\usepackage{amsfonts}
\usepackage{amsmath,amsthm}

\usepackage[all]{xy}
\usepackage{graphicx}

\usepackage{amssymb}
\usepackage{latexsym}

\DeclareMathAlphabet\EuFrak{U}{euf}{m}{n}	
\SetMathAlphabet\EuFrak{bold}{U}{euf}{b}{n}	

\newcommand{\ra}{\rightarrow}

\newcommand{\hra}{\hookrightarrow}

\newcommand{\ul}{\underline}
\newcommand{\wa}{\widehat}

\newcommand{\sC}{{\it C*}-}
\newcommand{\bC}{{\mathbb C}}

\newcommand{\bT}{{\mathbb T}}

\newcommand{\bZ}{{\mathbb Z}}
\newcommand{\bM}{{\mathbb M}}
\newcommand{\bN}{{\mathbb N}}

\newcommand{\ud}{{{\mathbb U}(d)}}
\newcommand{\sud}{{{\mathbb {SU}}(d)}}

\newcommand{\eps}{\varepsilon}

\newcommand{\mA}{\mathcal A}
\newcommand{\mB}{\mathcal B}
\newcommand{\mC}{\mathcal C}

\newcommand{\mE}{\mathcal E}
\newcommand{\mF}{\mathcal F}
\newcommand{\mG}{\mathcal G}
\newcommand{\mH}{\mathcal H}

\newcommand{\mK}{\mathcal K}
\newcommand{\mL}{\mathcal L}
\newcommand{\mM}{\mathcal M}
\newcommand{\mN}{\mathcal N}
\newcommand{\mO}{\mathcal O}
\newcommand{\mP}{\mathcal P}
\newcommand{\mQ}{\mathcal Q}

\newcommand{\mT}{\mathcal T}

\newcommand{\mU}{\mathcal U}
\newcommand{\mV}{\mathcal V}

\newcommand{\ii}{\iota,\iota}

\newcommand{\rr}{\rho,\rho}

\newcommand{\rs}{\rho,\sigma}
\newcommand{\rt}{\rho,\tau}
\newcommand{\sr}{\sigma,\rho}
\newcommand{\st}{\sigma,\tau}
\newcommand{\sss}{\sigma,\sigma}

\newcommand{\hrs}{H_d^r,H_d^s}

\newcommand{\tens}{ {\bf tens} }
\newcommand{\obc}{ {\bf obj} \ \mC }
\newcommand{\obt}{ {\bf obj} \ \mT }

\newtheorem{thm}{Theorem}[section]

\newtheorem{prop}[thm]{Proposition}

\newtheorem{defn}[thm]{Definition}

\theoremstyle{definition}
\newtheorem{ex}{Example}[section]

\theoremstyle{remark}

\numberwithin{equation}{section}

\begin{document}

\author{{\sf Ezio Vasselli}
                         \\{\it Dipartimento di Matematica}
                         \\{\it Universit\`a La Sapienza di Roma}
			 \\{\it P.le Aldo Moro, 2 - 00185 Roma - Italy }
 			 \\{\it (c/o Sergio Doplicher) }
                         \\{\sf vasselli@mat.uniroma2.it}}

\title{Group bundle duality,\\
       invariants for certain \sC algebras, \\
       and twisted equivariant $K$-theory}
\maketitle

\begin{abstract}
A duality is discussed for Lie group bundles vs. certain tensor \sC categories with non-simple identity, in the setting of Nistor-Troitsky gauge-equivariant $K$-theory. As an application, we study \sC algebra bundles with fibre a fixed-point algebra of the Cuntz algebra: a classification is given, and a cohomological invariant is assigned, representing the obstruction to perform an embedding into a continuous bundle of Cuntz algebras. Finally, we introduce the notion of twisted equivariant $K$-theory.

\bigskip

\noindent {\bf AMS Subj. Class.:} 19L47; 46L05; 46M15.


\end{abstract}

\section{Introduction.}
\label{intro}

An important branch of abstract harmonic analysis is group duality, whose task is to establish (possibly one-to-one) correspondences between groups and suitable tensor categories. For the case of a compact group, say $G$, the basic idea is to consider the tensor category of finite-dimensional, unitary representations of $G$; viceversa, given a tensor category $\mT$ (carrying some additional structure), one has to construct an embedding functor $F : \mT \hra {\bf hilb}$ into the category of finite-dimensional Hilbert spaces, and define $G$ as the group of natural transformations of $F$.

The existence of $F$ is a crucial step. In fact, in general $\mT$ may be not presented as a subcategory of the one of Hilbert spaces. For example, this happens in the case of algebraic quantum field theory, where the involved category is the one of superselection sectors of a \sC algebra of localized quantum observables. 
The duality for compact groups proved by Doplicher and Roberts in \cite{DR89} was motivated by such physical applications, and is the starting point for the work presented in the present paper.

The main feature of a tensor category $\mT$ is given by the fact that objects $\rs$ of $\mT$ can be multiplied via the tensor product $\rs \mapsto \rho \sigma$: the existence of a unit object $\iota$ is postulated, in such a way that $\rho \iota = \iota \rho = \rho$ for every object $\rho$. 
One of the properties assumed for the duality proved in \cite{DR89} is that $\iota$ is simple, i.e. that the algebra of arrows $(\ii)$ reduces to the complex numbers; in general, the tensor categories considered in the above-cited reference are such that $(\ii)$ is a commutative \sC algebra. 
This is an important point w.r.t. duality, because in such a case we have to look for an embedding functor into the category ${\bf mod} (\ii)$ of $(\ii)$-bimodules, instead of the one of Hilbert spaces.

In the present work, we study a class of tensor categories with non-simple unit object: we give a classification, and discuss the existence of an embedding functor into ${\bf mod} (\ii)$. In particular, we show that such embedding functors may not exist, in contrast to the case $(\ii) \simeq \bC$.

The Serre-Swan equivalence implies that we may equivalently consider the category of vector bundles over the spectrum of $(\ii)$, instead of ${\bf mod}(\ii)$. Since structural properties of the categories that we consider are encoded by geometrical invariants, we adopt this last point of view.

Our motivation arises from possible physical applications in low dimensional quantum field theory. Anyway, at a purely mathematical level we note that actions in the setting of vector bundles have been object of interesting research (\cite{NT04}).

The present paper is intended to be more self-contained as possible, and most of the material exposed here is a friendly exposition of \cite{Vas05}. The only original part is Sec.\ref{sec_tekt}, where we introduce a twisted equivariant topological $K$-functor.

Our work is organized as follows. In Sec.\ref{sec_kt}, we associate a $K$-theory group to a given \sC category; our construction provides usual \sC algebra $K$-theory in the case in which the given \sC category is the one of finitely generated Hilbert modules over a fixed \sC algebra. In Sec.\ref{sec_lie}, we expose a simple cohomological construction for compact Lie groups, which will be applied in the following sections. In Sec.\ref{sec_cbun}, we expose a version of the Serre-Swan equivalence for bundles of \sC algebras. In Sec.\ref{sec_tens}, we give a procedure to construct tensor categories with non-simple units, in the way as follows. Let $d \in \bN$, $\ud$ denote the unitary group, and $G \subseteq \ud$ a compact group. We denote by $NG \subseteq \ud$ the normalizer of $G$ in $\ud$, and by $QG := NG/G$ the quotient group. If $X$ is a compact Hausdorff space, we associate to every principal $QG$-bundle $\mQ$ a tensor category $T \mQ$ with unit $\iota$, such that $(\ii) \simeq C(X)$; roughly speaking, $T \mQ$ is a bundles of tensor categories over $X$, with fibre the category of tensor powers of the defining representation $G \hra \ud$. Let $G_{ab}$ denote the image of $G$ in the abelianized of $NG$; we define a cohomology class $\delta (\mQ) \in H^2(X,G_{ab})$, encoding the obstruction to construct an embedding functor $F : T \mQ \hra {\bf mod}(\ii)$. If $F$ does exist, then there is a $G$-bundle $\mG \ra X$ such that $T \mQ$ is the category of tensor powers of a $\mG$-equivariant vector bundle in the sense of \cite{NT04}. Existence of $F$ is equivalent to the condition that there is a principal $NG$-bundle, say $\mN$, such that $\mQ =$ $\mN \mod G$. In Sec.\ref{sec_b_og}, we provide a translation of the results of Sec.\ref{sec_tens} in terms of \sC algebra bundles. We consider the Cuntz algebra $\mO_d$, and the fixed point algebra $\mO_G$ w.r.t. a natural $G$-action; the invariant $\delta (\mQ)$ is now interpreted as the cohomological obstruction to embed a given $\mO_G$-bundle into some $\mO_d$-bundle.
In Sec.\ref{sec_tekt}, we associate to $\mQ$ a $K$-theory group $K_\mQ^0(X)$. If $\mQ = \mN \mod G$ for some principal $NG$-bundle $\mN$, then there exists a $G$-bundle $\mG \ra X$ such that $K_\mQ^0(X)$ can be embedded into the gauge equivariant $K$-theory $K_\mG^0(X)$ defined in the sense of \cite{NT04}.

{\bf Notation and background references.} Let $X$ be a compact Hausdorff space. We denote by $C(X)$ the \sC algebra of bounded, continuous, $\bC$-valued functions. If $x \in X$, we denote by $C_x(X)$ the closed ideal of functions vanishing on $x$. If $\left\{ X_i \right\}$ is an open cover of $X$, we denote $X_{ij} := X_i \cap X_j$, $X_{ijk} := X_{ij} \cap X_k$. In the present work, we adopt the convention that the set $\bN$ of natural numbers includes $0$. If $A$ is any set, then $id_A$ denotes the identity map on $A$. 
For topological $K$-theory, we refere to \cite{Kar}, while for (principal) group bundles our reference is \cite{Hus}. About \sC algebras and Hilbert modules, we refere to \cite{Bla}.

\section{\sC categories and $K$-theory.}
\label{sec_kt}

A \sC category $\mC$ is the data of a collection of objects, say $\obc$, such that for every pair $\rs \in \obc$ there is an associated Banach space $(\rs)$, called the space of arrows. A bilinear composition is defined, $(\st) \times (\rs)$ $\ra$ $(\rt)$, $t' , t \mapsto t' \circ t$, in such a way that $\left\| t' \circ t \right\|$ $\leq$ $\left\| t' \right\| \left\| t \right\|$, and an involution $* : (\rs)$ $\ra$ $(\sr)$ is assigned, in such a way that the \sC identity $\left\| t^* \circ t \right\| = \left\| t \right\|^2$ is satisfied. For every $\rho \in \obc$, we assume the existence of an identity $1_\rho \in (\rr)$ such that $1_\rho \circ t = t$, $t' \circ 1_\rho = t'$, $t \in (\sr)$, $t' \in (\rt)$. It is customary to postulate the existence of a {\em zero object} ${\bf o} \in$ $\obc$, such that $({\bf o} , \rho) =$ $\left\{ 0 \right\}$, $\rho \in \obc$. For basic properties about \sC categories, we refere the reader to \cite[\S 1]{DR89}.

As an immediate consequence of the above definition, we find that every $(\rr)$ is a unital \sC algebra w.r.t. the Banach space structure, composition of arrows and involution; moreover, every $(\rs)$ is a Hilbert $(\sss)$-$(\rr)$-bimodule w.r.t. composition of arrows and the $(\rr)$-valued scalar product $t,t' \mapsto t^* \circ t'$, $t,t' \in (\rs)$.
If $\rs \in \obc$, we say that $\rho$ is {\em unitarily equivalent} to $\sigma$ if there exists $u \in (\rs)$ with $u \circ u^* = 1_\sigma$, $u^* \circ u = 1_\rho$. We denote by $\left\{ \rho \right\} \subseteq \obc$ the set of objects which are unitarily equivalent to $\rho$.

Let $\mC$, $\mC'$ be \sC categories. A {\em functor} $\phi : \mC \ra \mC'$ is the data of a map $F : \obc \ra \obc'$ and a family $\left\{ F_{\rs} : (\rs) \ra ( F(\rho) , F(\sigma) ) \right\}$ of bounded linear maps preserving composition and involution. This implies that every $F_{\rr}$, $\rho \in \obc$, is a \sC algebra morphism.

\

We give some examples. A \sC category with a single object is clearly a unital \sC algebra. If $\mA$ is a unital \sC algebra, then the category with objects right Hilbert $\mA$-modules and arrows adjointable, right $\mA$-module operators is a \sC category, say ${\bf mod} (\mA)$. If $X$ is a compact Hausdorff space, the category ${\bf vect}(X)$ with objects vector bundles over $X$ and arrows morphisms of vector bundles is a \sC category (actually, by the Serre-Swan theorem \cite[Thm.I.6.18]{Kar} we may identify ${\bf vect}(X)$ with ${\bf mod}(C(X))$).

\

We now define some further structure. We say that $\mC$ has {\em subobjects} if for every $\rho \in \obc$ and projection $E = E^2 = E^* \in (\rr)$, there is $\sigma \in \obc$ and $S \in (\rs)$ such that $1_\sigma = S \circ S^*$, $E = S^* \circ S$. 
%
%
For every \sC category $\mC$, we may construct a 'larger' \sC category $\mC_s$ with subobjects, by defining ${\bf obj} \ \mC_s :=$ $\left\{ E = E^2 = E^* \in (\rr) , \rho \in \obc \right\}$, $(E,F)$ $:=$ $\left\{ t \in (\rs) : t = t \circ E = F \circ t \right\}$. $\mC_s$ is called the {\em closure for subobjects of} $\mC$.

$\mC$ has {\em direct sums}  if for every $\rs \in \obc$ there exists $\tau \in \obc$ with $\psi_\rho \in (\rt)$, $\psi_\sigma \in ( \st )$ such that $\psi_\rho \circ \psi_\rho^* + \psi_\sigma \circ \psi_\sigma^* = 1_\tau$, $\psi_\rho^* \circ \psi_\rho = 1_\rho$, $\psi_\sigma^* \circ \psi_\sigma = 1_\sigma$. The object $\tau$ is called {\em a direct sum} of $\rho$ and $\sigma$, and is unique up to unitary equivalence.
If a \sC category $\mC$ does not have direct sums, we can construct the {\em additive completition} $\mC_+$, having objects $n$-ples $\ul \rho := ( \rho_1 , \ldots , \rho_n )$, $n \in \bN$, and arrows spaces of matrices
$
(\ul \rho , \ul \sigma)
:= 
\left\{ (t_{ij}) : t_{ij} \in (\rho_j , \sigma_i ) \right\}
$.
\noindent In the case in which $\mC$ has direct sums, we define $\mC_+ := \mC$.
Note that the operations $\mC \mapsto \mC_s$ and $\mC \mapsto \mC_+$ do not commute, so that $\mC_{+,s}$ is not isomorphic to $\mC_{s,+}$: in general, there is an immersion
\[
\mC_{s,+} \hra \mC_{+,s} \ \ ,
\]
\noindent in fact every formal direct sum $E := ( E_1 , \ldots , E_n )$, $E_k \in {\bf obj} \ \mC_s$, $k = 1 , \ldots , n$, is a projection in the additive completition $\mC_+$.

\

Let $\mC$ be a \sC category with direct sums. Then, a semigroup $S(\mC)$ is associated with $\mC$, in the following way: for every $\rs \in \mC$, we consider a direct sum $\tau$ and the classes $\left\{ \rho \right\}$, $\left\{ \sigma \right\}$; then, we define
$
\left\{ \rho \right\} + \left\{ \sigma \right\} := \left\{ \tau \right\}
$.
Since $\left\{ \tau \right\}$ does not depend on the order of $\rs$, we find that $( S(\mC) , + )$ is an abelian semigroup with identity $\left\{ {\bf o} \right\}$.

\begin{defn}
\label{def_k}
Let $\mC$ be a \sC category. The $K$-theory group of $\mC$ is defined as the Grothendieck group, say $K_0(\mC)$, associated with $S(\mC_+)$.
\end{defn}

We give a class of examples: if $\mA$ is a unital \sC algebra, then $K_0({\bf mod}(\mA))$ coincides with the usual $K$-theory group $K_0(\mA)$. In particular, if $\mA = C(X)$, then $K_0({\bf mod}(C(X)))$ coincides with the topological $K$-theory $K^0(X)$. 

%
%

\section{A cohomological construction for certain Lie groups.}
\label{sec_lie}

Let $L$ be a compact group. A $L$-cocycle over $X$ is the data of a pair $\mL :=$ $(\left\{ X_i \right\},\left\{ g_{ij} \right\})$, where $\left\{ X_i \right\}$ is a finite open cover and $\left\{ g_{ij} \right\}$ is a family of continuous maps $g_{ij} : X_{ij} \ra L$ such that $g_{ik} (x) = g_{ij} (x) g_{jk}(x)$, $x \in X_{ijk}$. $L$-cocycles $(\left\{ X_i \right\},\left\{ g_{ij} \right\})$, $(\left\{ X'_l \right\},\left\{ g'_{lm} \right\})$ are said {\em equivalent} if there are continuous maps $u_{il} : X_i \cap X'_l \ra L$ such that $g_{ij} (x) = u_{il} (x) g'_{lm}(x) u_{jm}(x)^{-1}$, $x \in X_{ij} \cap X'_{lm}$. The set of equivalence classes of $L$-cocycles is called {\em cohomology set}, and is denoted by $H^1(X,L)$. It is well-known that $H^1(X,L)$ classifies principal $L$-bundles over $X$ (\cite[Chp.4]{Hus}); if $L$ is abelian, then $H^1(X,L)$ has a well-defined group structure, and coincides with the first cohomology group with coefficients in the sheaf of germs of continuous $L$-valued maps.
$H^1(X,\ \cdot\ )$ satisfies natural functoriality properties: if $\phi : L \ra L'$ is a group morphism, then the pair $( \left\{ X_i \right\} , \left\{ \phi \circ g_{ij} \right\} )$ define a $L'$-cocycle; thus, a map
\begin{equation}
\label{def_mc}
\phi_* : H^1(X,L) \ra H^1(X,L')
\end{equation}
\noindent is defined. If $L,L'$ are abelian, then $\phi_*$ is a group morphism.
As an example, let us consider the quotient $L_{ab}$ of $L$ w.r.t. the adjoint action $L \ra {\bf aut}L$. $L_{ab}$ is an abelian group, and there is a natural epimorphism $\pi_{L} : L \ra L_{ab}$, inducing a map $\pi_{L,*} : H^1(X,L) \ra H^1(X,L_{ab})$.

\

The following construction appeared in \cite[\S 4]{Vas05}. Let $d \in \bN$, $\ud$ denote the unitary group, and $G \subseteq \ud$ a compact group. We define $NG$ as the normalizer of $G$ in $\ud$, and $QG := NG / G$. Both $NG$ and $QG$ are compact Lie group; there is an epimorphism $p : NG \ra QG$ with kernel $G$, and a monomorphism $i_{NG} : G \hra NG$.
In general, the induced map 
\begin{equation}
\label{def_p*}
p_* : H^1(X,NG) \ra H^1(X,QG)
\end{equation}
\noindent is not surjective. We now define a cohomological class measuring the obstruction for surjectivity of $p_*$. Elementary computations (see \cite[\S 4]{Vas05} for details) show that there is a commutative diagram
\begin{equation}
\label{eq_cd_ab}
\xymatrix{
   1
   \ar[r]
 & G
   \ar[r]^-{ i }
   \ar[d]_-{ \pi_G }
 & NG
   \ar[r]^-{ p }
   \ar[d]^-{ \pi_{NG} }
 & QG
   \ar[r]
   \ar[d]^-{ \pi_{QG} }
 & 1
\\ 1
   \ar[r]
 & G_{ab}
   \ar[r]^-{ i_{ab} }
 & NG_{ab}
   \ar[r]^-{ p_{ab} }
 & QG_{ab}
   \ar[r]
 & 1
}
\end{equation}
\noindent where $1$ denotes the trivial group. By functoriality of $H^1(X,\ \cdot \ )$, and by applying the long exact sequence in sheaf cohomology, we obtain the commutative diagram
\begin{equation}
\label{eq_es}
\xymatrix{
   H^1 ( X , NG )
   \ar[r]^-{ p_* }
   \ar[d]^-{ \pi_{NG,*} }
 & H^1 ( X , QG )
   \ar[d]^-{ \pi_{QG,*} }
 & {}
\\ H^1 ( X , NG_{ab} )
   \ar[r]^-{ p_{ab,*} }
 & H^1 ( X , QG_{ab} )
   \ar[r]^-{ \delta_{ab} } 
 & H^2 ( X , G_{ab} )
}
\end{equation}

\begin{defn}
\label{def_ddc}
Let $d \in \bN$, and $G \subseteq \ud$ be a compact (Lie) group.
For every principal $QG$-bundle $\mQ \in H^1(X,QG)$, the Dixmier-Douady class of $\mQ$ is defined by
\[
\delta (\mQ) := \delta_{ab} \circ \pi_{QG,*} 
\ \in \ H^2 ( X , G_{ab} )
\ \ .
\]
\noindent If $\mQ$ belongs to the image of the map $p_* : H^1( X , NG ) \ra H^1( X , QG )$, then $\delta (\mQ) = 0$.
\end{defn}

We conclude the present section with some notation. $NG$ acts on $G$ {\it via} the adjoint action, and is naturally embedded in $\ud$. Thus, we have maps
\begin{equation}
\label{def_maps_g}
{\bf ad}  : NG \ra {\bf aut} G
\ \ , \ \
i_{\ud}   : NG \hra \ud
\end{equation}
\noindent which induce maps
\begin{equation}
\label{def_maps_h1g}
{\bf ad}_* : H^1(X,NG) \ra H^1(X,{\bf aut}G)
\ \ , \ \
i_{\ud,*}  : H^1(X,NG) \ra H^1(X,\ud) \ \ .
\end{equation}

\section{$C(X)$-algebras and \sC bundles.}
\label{sec_cbun}

In the present section, we expose some basic properties of \sC algebra bundles. Such properties will be applied in Sec.\ref{sec_b_og}.

Let $X$ be a compact Hausdorff space.
A unital \sC algebra $\mA$ is said $C(X)${\em -algebra} if there is a unital morphism $C(X) \ra \mA' \cap \mA$; in the sequel, elements of $C(X)$ will be identified with their image in $\mA$, so that $C(X)$ may be regarded as a unital \sC subalgebra of $\mA' \cap \mA$. For every $x \in X$, we consider the ideal $I_x := \left\{ fa : f \in C_x(X) , a \in \mA \right\}$, and define the {\em fibre} $\mA_x := \mA / I_x$ with the epimorphism $\pi_x : \mA \ra \mA_x$; this allows to regard at every $a \in \mA$ as a vector field $\wa a :=$ $(\pi_x (a))_{x \in X} \in$ $\prod_x \mA_x$, $a \in \mA$. In general, the {\em norm function} $n_a (x) := \left\| \pi_x(a) \right\|$ , $x \in X$, is upper-semicontinuous. In the case in which every $n_a$, $a \in \mA$, is continuous, then $\mA$ is called {\em continuous bundle} over $X$. A $C(X)${\em -algebra morphism} $\phi : \mA \ra \mB$ is a \sC algebra morphism such that $\phi (fa) = f\phi(a)$, $f \in C(X)$, $a \in \mA$. If $\mF$ is a fixed \sC algebra, we denote by ${\bf bun}(X,\mF)$ the set of $C(X)$-isomorphism classes of $C(X)$-algebras having fibres isomorphic to $\mF$. 
Good references about $C(X)$-algebras are \cite{Bla96,Nil96}.

There is categorical equivalence between $C(X)$-algebras and a class of topological objects, called \sC bundles (see \cite[Thm.5.13]{Gie82}, \cite[\S 10.18]{Gie82}). A {\em \sC-bundle} is the data of a Hausdorff space $\Sigma$ endowed with a surjective, open, continuous map $Q : \Sigma \ra X$ such that every fibre $\Sigma_x := Q^{-1}(x)$, $x \in X$, is homeomorphic to a unital \sC algebra with identity $1_x$. $\Sigma$ is required to be full, i.e. for every $\sigma \in \Sigma$ there exists a continuous section $a : X \hra \Sigma$, $a \circ Q = id_X$, such that $a \circ Q (\sigma) = \sigma$. 
Let $Q' : \Sigma' \ra X$ be a \sC bundle. A \sC {\em bundle morphism} from $\Sigma$ into $\Sigma'$ is the data of a continuous map $\phi : \Sigma \ra \Sigma'$ such that 
\begin{enumerate}
\item  $Q' \circ \phi = Q$; this implies that $\phi (\Sigma_x) \subseteq \Sigma'_x$, $x \in X$;
\item  $\phi_x := \phi |_{\Sigma_x} : \Sigma_x \ra \Sigma'_x$ is a \sC algebra morphism for every $x \in X$.

\end{enumerate}

\noindent $\phi$ is said {\bf isomorphism} if every $\phi_x$, $x \in X$, is a \sC algebra isomorphism.

The set $S_X(\Sigma)$ of continuous sections of $\Sigma$ is endowed with a natural structure of $C(X)$-algebra: for every $a,a' \in \mA$, we define $a+a'$, $aa'$, $a^*$ as the maps $(a+a') (x) :=$ $a(x) + a'(x)$, $(aa')(x) :=$ $a(x) a'(x)$, $a^*(x) := a(x)^*$, $x \in X$, and introduce the norm $\left\| a \right\| :=$ $\sup_x \left\| a(x) \right\|$. Moreover, for every $f \in C(X)$, it turns out that the map $f1(x) := f(x) 1_x$, $x \in X$, belongs to $S_X(\Sigma)$.

Viceversa, if $\mA$ is a unital $C(X)$-algebra, then the set $\wa \mA :=$ $\bigsqcup_{x \in X} \mA_x$ is endowed with a natural surjective map $Q : \wa \mA \ra X$, $Q \circ \pi_x (a) := x$. For every open $U \subseteq X$, $a \in \mA$, $\eps > 0$, we consider the following subset of $\wa \mA$:
\begin{equation}
\label{def_tube}
T_{U,a,\eps} :=
\left\{ 
\sigma \in \wa \mA 
\ : \ 
Q(\sigma) \in U
\ {\mathrm{and}} \
\left\| \sigma - \pi_{Q(\sigma)} (a) \right\| < \eps
\right\} \ .
\end{equation}
\noindent The family $\left\{ T_{U,a,\eps} \right\}$ provides a basis for a topology on $\wa \mA$. It can be proved that $Q : \wa \mA \ra X$ is a \sC bundle, and that there is an isomorphism $\mA \ra S_X(\wa \mA)$, $a \mapsto \wa a :=$ $\left\{ X \ni x \mapsto \pi_x(a) \right\}$.

\begin{thm}
Let $X$ be a compact Hausdorff space. The map $\mA \mapsto \wa \mA$ provides an equivalence between the category of $C(X)$-algebras (with arrows $C(X)$-algebra morphisms) and the category of \sC bundles over $X$ (with arrows \sC bundle morphisms).
\end{thm}

Some examples follow. If $\mF$ is a unital \sC algebra, then $\mA := C(X) \otimes \mF$ is a $C(X)$-algebra with \sC bundle $\wa \mA = X \times \mF$, endowed with the projection  $Q(x,b) := x$, $x \in X$, $b \in \mF$; $\wa \mA$ is called the {\em trivial} \sC bundle. Let $C(X) \hra C(Y) =: \mA_Y$ be an inclusion of unital, abelian \sC algebras; then, a surjective map $q : Y \ra X$ is defined, and the \sC bundle $\wa \mA_Y \ra X$ has fibres $\mA_x = C(q^{-1}(x))$, $x \in X$. Let $d \in \bN$, and $\mE \ra X$ denote a rank $d$ vector bundle; then, the \sC algebra $L(\mE)$ of endomorphisms of $\mE$ is a continuous bundle of \sC algebras over $X$, with fibres isomorphic to the matrix algebra $\bM_d$.

\

We introduce some notation and terminology. Let $L$ denote a compact group acting by automorphisms on $\mF$. Moreover, let $\mA$, $\mA'$ be $C(X)$-algebras with fibre $\mF$ (i.e., $\mA_x \simeq$ $\mA'_x \simeq$ $\mF$, $x \in X$), and $\beta : \mA \ra \mA'$ a $C(X)$-algebra isomorphism. Then, for every $x \in X$ an automorphism $\beta_x \in {\bf aut}\mF$ is defined, in such a way that $\beta_x \circ \pi_x (a) =$ $\pi'_x \circ \beta (a)$, $a \in \mA$. If $\beta_x$ belongs to the image of the $L$-action for every $x \in X$, then we say that $\beta$ is $L${\em -covariant}, and use the notation
\begin{equation}
\label{eq_geq}
\beta :  \mA  \ra_L \mA' \ .
\end{equation}
\noindent Let now $I : \mF_0 \hra \mF$ be an inclusion of unital \sC algebras. Suppose that there are $C(X)$-algebras $\mA_0 \in$ ${\bf bun}(X,\mF_0)$, $\mA \in$ ${\bf bun}(X,\mF)$, and a a $C(X)$- monomorphism $\phi : \mA_0 \ra \mA$ such that $\phi_x = I$ for every $x \in X$. Then, we say that $\phi$ is $I$-{\em covariant}, and use the notation
\begin{equation}
\label{eq_ieq}
\phi :  \mA_0  \hra_{I} \mA \ \ .
\end{equation}

We conclude the present section by presenting a construction for continuous bundles having as fibre a fixed \sC dynamical system. 
Let $L$ be a compact group, and $\mF$ a unital \sC algebra with an automorphic action $\alpha : L \ra {\bf aut} \mF$. For every $L$-cocycle $\mL := $ $(\left\{ X_i \right\},\left\{ g_{ij} \right\})$, we define a \sC bundle $Q : \wa \mA_\mL \ra X$ as the {\em clutching} of the family of trivial bundles $\left\{ X_i \times \mF \right\}$ w.r.t. the maps $\left\{ \alpha \circ g_{ij} : \right.$ $X_{ij} \ra$ $\left. {\bf aut} \mF \right\}$ (in the same way as in \cite[I.3.2]{Kar}). We denote by $\mA_{\mL}$ the $C(X)$-algebra of continuous sections associated with $\wa \mA_\mL$.
If $\mH$ is a $L$-cocycle equivalent to $\mL$, then it is easily verified that there is an isomorphism $\wa \mA_{\mL} \simeq \wa \mA_{\mH}$. In such a way, we defined a map 
\begin{equation}
\label{def_cb}
\alpha_* : H^1(X,L) \ra {\bf bun}(X,\mF)
\ , \
\mL \mapsto \mA_{\mL} \ .
\end{equation}
\noindent In general $\alpha_*$ is not injective, unless $L \simeq {\bf aut} \mF$; in fact, cocycle equivalence gives rise to an $L$-covariant isomorphism.

\section{Tensor \sC categories.}
\label{sec_tens}

A {\em tensor} \sC {\em category} is a \sC category $\mT$ endowed with a bifunctor $\otimes : \mT \times \mT \ra \mT$, called the {\em tensor product}. In explicit terms, for every pair $\rs \in \obt$ there is an object $\rho \sigma \in \obt$; for every $\rho' , \sigma' \in \obt$, there are bilinear maps $(\rho , \sigma) \times (\rho' , \sigma')$ $\ra$ $( \rho \rho' , \sigma \sigma' )$, $t,t' \mapsto t \otimes t'$. The existence of an {\em identity object} $\iota \in \obt$ is postulated, in such a way that $\iota \rho = \rho \iota = \rho$, $\rho \in \obt$, $t = t \otimes 1_\iota = 1_\iota \otimes t$, $t \in (\rs)$. The data of a tensor \sC category with identity object $\iota$ will be denoted by the triple $( \mT , \otimes , \iota )$. For basic notions on tensor \sC categories, we refere to \cite[\S 1]{DR89}.

It is a consequence of the above definition that the \sC algebra $(\ii)$ is abelian; we will denote by $X^\iota$, the (compact, Hausdorff) spectrum of $(\ii)$, so that there is an identification $(\ii) \simeq C(X^\iota)$.

\

Well-known examples of tensor \sC categories are the one of Hilbert spaces, say ${\bf hilb}$, endowed with the usual tensor product, and the one of vector bundles ${\bf vect}(X)$ over a compact Hausdorff space $X$. In the first case $X^\iota$ reduces to a single point ($\iota = \bC$, so that $(\ii) \simeq \bC$); in the second case $X^\iota = X$ ($\iota = X \times \bC$).

\

The eventual commutativity up-to-unitary-equivalence of the tensor product is described by the property of {\em symmetry}. A tensor \sC category $(\mT , \otimes , \iota )$ is said {\em symmetric} if for each $\rs \in \obt$ there is a unitary 'flip' $\eps_{\rho \sigma} \in (\rho \sigma , \sigma \rho)$ such that $\eps_{\sigma \sigma'} \circ (t \otimes t') =$ $(t' \otimes t) \circ \eps_{\rho \rho'}$, $t \in (\rs)$, $t' \in (\rho' , \sigma')$. A symmetric tensor \sC category is denoted by $( \mT , \otimes , \iota , \eps )$.

The tensor \sC categories ${\bf hilb}$, ${\bf vect}(X)$ are symmetric. Another basic example is given by the {\em dual} of a compact group $G$, i.e. the category with objects unitary, finite-dimensional representations of $G$.

\subsection{Duals of compact Lie groups.}

Let $G$ be a compact group endowed with a faithful representation over a rank $d$ Hilbert space $H_d$, $d \in \bN$. We regard at $G$ as a compact Lie subgroup of the unitary group $\ud$.
We define $\wa G$ as the tensor \sC category with objects the $r$-fold tensor powers $H_d^r$, $r \in \bN$ (for $r = 0$ we define $\iota :=$ $H_d^0 :=$ $\bC$), and arrows the spaces $( \hrs )_G$ of linear operators $t :$ $H_d^r \ra$ $H_d^s$ such that $g_s \circ t \circ g_r^* = t$ $\forall g \in G$, where
\begin{equation}
\label{def_gr}
g_r := \otimes^r g \ \in \ H_d^r \ .
\end{equation}

It is well-known that by performing the closure for subobjects and the additive completition of $\wa G$, we obtain all the finite-dimensional representations of $G$. $\wa G$ is symmetric, in fact it is endowed with the flip operators $\theta_{r,s} \in$ $(H_d^{r+s},H_d^{r+s})_G$, $\theta_{r,s}(\psi \otimes \psi') := \psi' \otimes \psi$, $\psi \in H_d^r$, $\psi \in H_d^s$. Thus, we have a symmetric tensor \sC category $( \wa G , \otimes , \iota , \theta  )$.

A {\em symmetric autofunctor} of $\wa G$ is the data of a family $F$ of Banach space isomorphisms
\[
F^{r,s}  : (\hrs)_G \ra (\hrs)_G  \ \ , \ \ r,s \in \bN
\]
\noindent such that $F^{r,s}(t \circ t') =$ $F^{l,s}(t) \circ F^{r,l}(t')$, $F^{r,s}(t^*) =$ $F^{r,l}(t)^*$, $F^{r+r',s+s'}(t \otimes t'') =$ $F^{r,s}(t) \otimes F^{r',s'}(t'')$, $F^{r+s,r+s}(\theta_{r,s}) = \theta_{r,s}$, $t \in (\hrs)_G$, $t' \in ( H_d^{r'} , H_d^{s'} )_G$, $t'' \in ( H_d^l , H_d^s )_G$. The set ${\bf aut}_\theta \wa G$ of symmetric autofunctors of $\wa G$ is endowed with a group structure w.r.t. the composition $F,G \mapsto G \circ F :=$ $\left\{ G^{r,s} \circ F^{r,s}  \right\}$ and inverse $F^{-1} :=$ $\left\{ (F^{r,s})^{-1} \right\}$. 
 
Now, every $u \in NG$ defines maps 
\[
\wa u^{r,s} : (\hrs)_G \ra (\hrs)_G
\ \ , \ \ 
\wa u^{r,s} (t) := u_s \circ t \circ u_r^* 
\]
\noindent (the term $u_s$ is defined according to (\ref{def_gr})). A direct check show that that the family $\wa u :=$ $\left\{ \wa u^{r,s} \right\}_{r,s}$ defines an element of ${\bf aut}_\theta \wa G$. In particular, we find
\begin{equation}
\label{eq_theta}
\wa g^{r+s,r+s}(\theta_{r,s}) = \theta_{r,s}
\end{equation}
\noindent (see \cite[\S 2]{DR87}). Since by definition 
\begin{equation}
\label{eq_tggt}
\wa g^{r,s}(t) = t 
\ \ , \ \
g \in G \ , \ t \in (\hrs)_G \ ,
\end{equation}
\noindent we conclude that
\begin{equation}
\label{eq_ug}
\wa u = \wa{ug} \ , \ u \in NG \ , \ g \in G \ \ ,
\end{equation}
\noindent so that $\left\{ u \mapsto \wa u  \right\}$ factorizes through a map
\begin{equation}
\label{eq_wg}
QG \ra {\bf aut}_\theta \wa G 
\ \ , \ \
y \mapsto \wa y \ .
\end{equation}

%
%

\subsection{Special categories and group bundles.}
\label{ssec_sc}

%
Let $X$ be a compact Hausdorff space, $d \in \bN$, $G \subseteq \ud$ a compact group. A {\em special category} is the data of a symmetric tensor \sC category $( \mT , \otimes , \iota , \eps )$ with objects the positive integers $r \in \bN$, and arrows the Banach $C(X)$-bimodules $\mM_{r,s}$ of continuous sections of vector bundles $\mE_{r,s} \ra X$ with fibre $(\hrs)_G$, $r,s \in \bN$ (the left $C(X)$-action of $\mM_{r,s}$ is assumed to coincide with the right one). The tensor product is defined as follows:
\[
\left\{
\begin{array}{ll}
r,s \mapsto r+s 
\ \ , \ \ 
r,s \in \bN
\\
t ,t' \mapsto t \otimes_X t'  \in  \mM_{r+r',s+s'}
\ \ , \ \ 
t \in \mM_{r,s} \ , \ t' \in \mM_{r',s'} \ ,
\end{array}
\right.
\]
\noindent where $\otimes_X$ denote the tensor product in the category of Banach $C(X)$-bimodules (\cite[Chp.VI]{Bla}). 
%
%
Note that $\iota := 0 \in \bN$ is the identity object, with $(\ii) = C(X)$. We denote by $\tens ( X , \wa G )$ the set of isomorphism classes of special categories having spaces of arrows with fibres $(\hrs)_G$, $r,s \in \bN$.

As an example, we consider the {\em trivial special category} $X \times \wa G$ with arrows $\mM_{0 , r,s} = C(X) \otimes (\hrs)_G$, $r,s \in \bN$, and symmetry $\theta^X :=$ $\left\{ \theta^X_{r,s} := \right.$ $\left. 1_X \otimes \theta_{r,s} \right\}$, where $1_X \in C(X)$ denotes the identity. Note that for every $r,s \in \bN$, it turns out that $\mM_{0,r,s}$ is the module of sections of $\mE_{0 , r,s} :=$ $X \times (\hrs)_G$.

\

We now give a simple procedure to construct special categories. Let $\mQ$ be a principal $QG$-bundle with associated cocycle $( \left\{ X_i \right\} , \left\{ y_{ij} \right\} )$. For every $r,s \in \bN$, we denote by ${\bf aut} (\hrs)_G$ the (topological) group of isometric linear maps of $(\hrs)_G$. By composing with the isomorphism (\ref{eq_wg}), we obtain maps
\[
\wa y^{r,s}_{ij} : X_{ij} \ra {\bf aut} (\hrs)_G
\ \ , \ \
x \mapsto \wa {y}_{ij} (x)^{r,s}
\ \ ,
\]
\noindent which define ${\bf aut} (\hrs)_G$-cocycles. We denote by $\mE_{ \mQ ,r,s} \ra X$ the vector bundle with fibre $(\hrs)_G$ and transition maps $\left\{ \wa y^{r,s}_{ij} \right\}$. Now, for every $r,s$, $r',s' \in \bN$, it turns out 
\[
\wa {y}_{ij} (x)^{r,s} \otimes \wa {y}_{ij} (x)^{r',s'}
=
\wa {y}_{ij} (x)^{r+r',s+s'} 
\ , \
x \in X_{ij}
\ ;
\]
\noindent this implies that there are inclusions
\begin{equation}
\label{eq_vbtp} 
\mE_{ \mQ ,r,s}  \otimes \mE_{ \mQ ,r',s'} \subseteq \mE_{ \mQ ,r+r',s+s'} \ \ ,
\end{equation}
\noindent where $\otimes$ stands for the tensor product in ${\bf vect} (X)$. We denote by $\mM_{ \mQ ,r,s}$ the module of continuous sections of $\mE_{ \mQ ,r,s}$, and define the following \sC category:
\[
T \mQ \ := \
\left\{
\begin{array}{ll}
{\bf obj} \ T \mQ := \bN
\\
(r,s)_\mQ := \mM_{ \mQ ,r,s} \ \ , \ \ r,s \in \bN \ \ .
\end{array}
\right.
\]
\noindent The relations (\ref{eq_vbtp}) imply that $\mM_{ \mQ ,r,s}  \otimes_X \mM_{ \mQ ,r',s'} \subseteq$ $\mM_{ \mQ ,r+r',s+s'}$, so that $T \mQ$ is a tensor \sC category. Note that $(0,0)_\mQ = C(X)$; more in general, for every $r \in \bN$ it turns out that $\mM_{\mQ,r,r}$ is a $C(X)$-algebra with fibre $(H_d^r,H_d^r)_G$. About the symmetry, let us consider the constant maps $\eps^i_{r,s} (x) := \theta_{r,s}$, $x \in X_i$; then, (\ref{eq_theta}) implies that
\[
\wa y^{r,s}_{ij} \ (\eps^i_{r,s}) |_{X_{ij}} = \eps^j_{r,s} |_{X_{ij}} \ .
\]
\noindent The previous relations imply that we can glue the local sections $\eps^i_{r,s} : X_i \ra$ $X_i \times (H_d^{r+s},H_d^{r+s})_G$ by using the transition maps $\wa y^{r,s}_{ij}$, and obtain elements $\eps_{r,s} \in$ $\mM_{\mQ,r+s,r+s}$. Some routine computations show that  the family $\left\{ \eps_{r,s} \right\}$ defines a symmetry for $T \mQ$. Thus, $T \mQ \in \tens (X,\wa G)$. It is easy to prove that if $\mQ' \in$ $H^1(X,QG)$ is cocycle equivalent to $\mQ$ then there is an \sC category isomorphism $T \mQ \simeq T \mQ'$ preserving symmetry and tensor product. Thus, we defined a map
\[
T : H^1 (X,QG) \ra \tens (X,\wa G) \ .
\]

\begin{ex}[Tensor powers of a vector bundle]
\label{ex_d_vb}
Let us consider the trivial group $G := {\mathbb I}_d$, so that $NG = QG = \ud$. Every $\ud$-cocycle $\mU :=$ $(\left\{ X_i \right\} , \left\{ u_{ij} \right\})$ can be regarded as the set of transition maps of a rank $d$ vector bundle $\mE \ra X$. For every $r \in \bN$, we denote by $\mE^r$ the $r$-fold tensor power of $\mE$. The tensor \sC category $T \mU$ has spaces of arrows the bimodules of continuous sections of vector bundles $\mE_{\mU , r , s} \ra X$, having fibre $(\hrs)$ and transition maps $\wa u_{ij}^{r,s}$. It is well-known that every $\mE_{\mU , r , s}$ can be identified with the vector bundle of morphisms from $\mE^r$ into $\mE^s$ (see \cite[I.4.8(c)]{Kar}). Thus, $T \mU$ is isomorphic to the full \sC subcategory of ${\bf vect}(X)$ with objects the tensor powers of $\mE$. We denote by $\wa \mE$ such tensor \sC category.
\end{ex}

\begin{ex}[The dual of a Lie group bundle]
\label{ex_d_gb}
We recall the reader to the notation of Sec.\ref{sec_lie} (in particular, the maps (\ref{def_maps_g},\ref{def_maps_h1g}) ). Let $\mN :=$ $( \left\{ X_i \right\} , \left\{ u_{ij} \right\} ) \in$ $H^1(X,NG)$. We consider the vector bundle $\mE \ra X$ with transition maps $\mU :=$ $i_{\ud,*} \mN \in$ $H^1(X,\ud)$, and the $QG$-cocycle $\mQ := $ $p_*(\mN) \in$ $H^1(X,QG)$. The symmetric tensor \sC category $T \mQ$ has arrows the bimodules of continuous sections of vector bundles $\mE_{\mQ,r,s} \ra X$, which have fibre $(\hrs)_G$ and transition maps $p \circ \wa{u}_{ij}^{r,s}$. Let now $\mE_{\mU,r,s} \ra X$ be the vector bundles associated with $\mU$, according to the previous example. Then, we may regard at each $\mE_{\mQ,r,s}$ as a vector subbundle of $\mE_{\mU,r,s}$; in fact, by (\ref{eq_ug}), every set of transition maps $\left\{ \wa u_{ij}^{r,s} \right\}$ acts in the same way as $\left\{ p \circ \wa{u}_{ij}^{r,s} \right\}$ on the trivial bundles $X_i \times (\hrs)_G \subseteq$ $X_i \times (\hrs)$.
Thus, $T \mQ$ is a symmetric tensor \sC subcategory of $\wa \mE$.
Let now $\mG \ra X$ be the group bundle with fibre $G$ and transition maps ${\bf ad}_* \mN$. Since ${\bf ad} (u) = \wa u^{1,1}$, $u \in NG$, there is an inclusion $\mG \subset \mE_{\mU,1,1}$. Since $\mE_{\mU,1,1}$ is the vector bundle of endomorphisms of $\mE$, we find that $\mE$ is $\mG$-equivariant in the sense of \cite[\S 1]{NT04}, i.e. there is an action
\begin{equation}
\label{eq_gba}
\mG \times_X \mE \ra \mE
\ \ , \ \ 
(g,v) \mapsto g(v) \ ,
\end{equation}
\noindent where $\mG \times_X \mE$ denotes the fibered cartesian product. Note that in (\ref{eq_gba}) the base space of $\mE$ coincides with the base space of $\mG$; in \cite{NT04}, the base space of $\mE$ is a topological bundle carrying a $\mG$-action.
Now, for every $r,s \in \bN$ there is an action
\[
\mG \times_X \mE_{\mU,r,s} \ra \mE_{\mU,r,s} 
\ \ , \ \ 
(g,t) \mapsto \wa g^{r,s} (t) \ .
\]
\noindent In particular, $t \in \mE_{\mQ,r,s}$ if and only if $\wa g^{r,s} (t) = t$ for every $g \in \mG$, in the same way as in (\ref{eq_tggt}); in such a case, we say that $t$ is $\mG${\em -equivariant}. We conclude that $T \mQ$ is the tensor \sC subcategory of $\wa \mE$ with arrows the bimodules of $\mG$-equivariant morphisms between tensor powers of $\mE$. We denote by $\wa \mE_\mG$ such tensor \sC category.
\end{ex}

\

Let $T \mQ$ be the special category associated with $\mQ \in H^1(X,QG)$. An {\em embedding functor} is the data of a rank $d$ vector bundle $\mE \ra X$ and a \sC monofunctor $F : T \mQ \hra \wa \mE$, preserving tensor product and symmetry. In the following results, we characterize tensor \sC categories $T \mQ$ admitting an embedding functor, and give a cohomological obstruction for the existence. This provides a duality between group bundles and special categories admitting an embedding functor.

\begin{thm}
[\cite{Vas05}, Thm.7.1, Thm.7.3]
\label{thm_class}
Let $G \subseteq \sud$. The following are equivalent:
\begin{enumerate}
\item  there exists an embedding functor $F : T \mQ \hra \wa \mE$, and a fibred $G$-bundle $\mG \ra X$ such that $\mE$ is $\mG$-equivariant and $T \mQ = \wa \mE_\mG$;
\item there exists a principal $NG$-bundle $\mN$ such that $p_* \mN = \mQ$.
\end{enumerate}
\end{thm}

The interplay between the above objects is the following: if there exists $\mN \in H^1(X,NG)$ with $p_* \mN = \mQ$, then $\mE$ is defined by the cocycle $i_{\ud,*} \mN \in$ $H^1(X,\ud)$, and $\mG$ is defined by the cocycle ${\bf ad}_* \mN \in$ $H^1(X,{\bf aut}G)$. The condition $G \subseteq \sud$ is motivated by the fact that the proof of the previous theorem lies on the notion of special object (\cite[\S 3, Lemma 6.7]{DR89}); it is our opinion that it should suffice to assume $G \subseteq \ud$, and this point is object of a work in progress. 
As a direct consequence of the previous theorem, we obtain
\begin{thm}
\label{thm_ci}
For every $T \mQ \in \tens (X , \wa G)$, we define
\[
\delta (T \mQ) := \delta (\mQ) \in H^2( X , G_{ab} )
\ \ 
( {\mathrm{ the \ Dixmier-Douady \ class \ of \ }} T \mQ) \ .
\]
\noindent If there exists an embedding functor $F : T \mQ \hra \wa \mE$, then $\delta(T \mQ) = 0$. In particular, if $T \mQ$ is the trivial special category (i.e., $\mQ \simeq X \times QG$) then $\delta (T \mQ) = 0$. 
\end{thm}

It is not difficult to construct special categories with $\delta (T \mQ) \neq 0$. For example, let us suppose that the epimorphism $\pi_{QG} :$ $QG \ra$ $QG_{ab}$ admits a left inverse $S :$ $QG_{ab} \hra$ $QG$, $\pi_{QG} \circ S = id_{QG_{ab}}$, and let $X$ be a space such that there is $z \in$ $H^1(X,QG_{ab})$ with $\delta_{ab}(z) \neq 0$. We define $\mQ :=$ $S_* z \in$ $H^1(X,QG)$; by construction, $\delta (T \mQ) =$ $\delta (\mQ) =$ $\delta_{ab} (z) \neq$ $0$. Explicit examples are given in \cite[\S 7.0.8]{Vas05}.

\section{\sC bundles with fibre $\mO_G$.}
\label{sec_b_og}


The following construction appeared in \cite{DR87}, and can be interpreted as a \sC algebraic version of the Tannaka duality.
Let $d \in \bN$; we denote by $\mO_d$ the Cuntz algebra (\cite{Cun77}) generated by a multiplet $\left\{ \psi_i \right\}_{i=1}^d$ of isometries satisfying the relations
\begin{equation}
\label{rel_cuntz}
\psi_i^* \psi_j = \delta_{ij} 1 
\ \ , \ \
\sum_i \psi_i \psi_i^* = 1
\ \ .
\end{equation}
\noindent Let us denote by $H_d \subset \mO_d$ the Hilbert space spanned by $\left\{ \psi_i \right\}$, endowed with the scalar product $\left \langle \psi , \psi' \right \rangle 1 :=$ $\psi^* \psi'$, $\psi , \psi \in H_d$; then, it is clear that $H_d$ is isomorphic to the canonical rank $d$ Hilbert space. Let now $I := \left\{ i_1 , \ldots , i_r \right\}$ be a multiindex with lenght $|I| := r \in \bN$; we introduce the notation
\[
\psi_I := \prod_{k=1}^r \psi_{i_k} \in \mO_d
\ \ , \ \ 
\]
\noindent and denote by $H_d^r$ the vector space spanned by $\left\{ \psi_I \right\}$, which we identify with the $r$-fold tensor power of $H_d$. With the above notation, the Banach space
\begin{equation}
\label{def_hrs}
(\hrs) := {\mathrm{span}} \ \left\{ 
\psi_I \psi_J^* \ , \ |I| = s , |J| = r
\right\}
\ , \
r,s \in \bN \ ,
\end{equation}
\noindent can be naturally identified with the set of linear operators from $H_d^r$ into $H_d^s$.
%
%
%
%
%
%
By identifying $H_d$ with the canonical rank $d$ Hilbert space, we get a natural $\ud$-action $\ud \times H_d \ra H_d$, $u,\psi \mapsto u \psi$. Universality of the Cuntz algebra implies that there is an automorphic action
\begin{equation}
\label{def_ua}
\ud \ra {\bf aut} \mO_d
\ \ , \ \ 
u \mapsto \wa u \ \ ,
\end{equation}
\noindent $\wa u (\psi) := u \psi$, $u \in \ud$, $\psi \in H_d$. If we restrict (\ref{def_ga}) to elements of a compact group $G \subseteq \ud$, then we get an action
\begin{equation}
\label{def_ga}
G \ra {\bf aut} \mO_d
\ \ , \ \ 
g \mapsto \wa g \ \ .
\end{equation}
\noindent We denote by $\mO_G$ the fixed-point algebra of $\mO_d$ w.r.t. the action (\ref{def_ga}), and by
\[
i_G : \mO_G \hra \mO_d
\]
\noindent the natural inclusion. Since $\wa g (t) = g_s \circ t \circ g_r^*$, $g \in G$, $t \in (\hrs)$, we may identify $(\hrs)_G$ with $(\hrs) \cap \mO_G$.

Let us now consider the normalizer $NG \subseteq \ud$. Then, every $u \in NG$ defines an automorphism $\wa u \in {\bf aut} \mO_G$: in fact, if $g \in G$, $t \in \mO_G$, then $\wa g ( \wa u (t)) =$ $\wa u \circ \wa{g'} (t) =$ $\wa u (t)$, where $g' \in G$. It is clear that $\wa u = \wa{ug}$ for every $g \in G$, thus there is an automorphic action
\begin{equation}
\label{def_qa}
QG \ra {\bf aut} \mO_G
\ \ , \ \
y \mapsto \wa y \ \ .
\end{equation}

%
%


We now apply the construction (\ref{def_cb}) to the actions (\ref{def_ua}), (\ref{def_qa}), and obtain maps
\begin{equation}
\label{def_bu}
H^1(X,\ud) \ra {\bf bun}(X,\mO_d)
\ \ , \ \
\mU \mapsto \mO_\mU \ \ ,
\end{equation}
\begin{equation}
\label{def_bq}
H^1(X,QG) \ra {\bf bun}(X,\mO_G)
\ \ , \ \
\mQ \mapsto \mO_\mQ \ \ .
\end{equation}
\noindent Let $\mE \ra X$ be the vector bundle with associated $\ud$-cocycle $\mU$, and $S \mE$ the Hilbert $C(X)$-bimodule of continuous sections of $\mE$. It is proved in \cite[Prop.4.2]{Vas} that $\mO_\mU$ is the Cuntz-Pimsner algebra associated with $S \mE$. The map (\ref{def_bu}) is not injective: if $X$ is a finite-dimensional $CW$-complex, in order to obtain an isomorphism $\mO_\mU \simeq \mO_\mV$ it suffices that $\mU , \mV \in$ $H^1(X,\ud)$ are transition maps of vector bundles having the same class in $K^0(X)$
(\cite[Prop.10]{Vas04}).

The following result is a translation of Thm.\ref{thm_class}, Thm.\ref{thm_ci} in terms of $C(X)$-algebras (note that we use the notation (\ref{eq_geq},\ref{eq_ieq})).

\begin{prop}
\label{prop_cog}{\bf \cite[Prop.7.11]{Vas05}}
$\mQ = \mQ'$ $\in H^1 ( X,QG )$ if and only if there is a $QG$-covariant $C(X)$-isomorphism $\alpha : \mO_\mQ  \ra_{QG}$ $\mO_{\mQ^\prime}$.
Moreover, the class 
\[
\delta (\mO_\mQ) := \delta(\mQ) \in H^2( X,G_{ab} )
\]
\noindent measures the obstruction to find a vector bundle $\mE \ra X$ with transition maps $\mU \in H^1(X,\ud)$ implementing a $C(X)$-monomorphism $\mO_\mQ \hra_{i_G} \mO_\mU$: if $\delta (\mO_\mQ) \neq 0$, then such a vector bundle does not exist.
\end{prop}

\section{Twisted equivariant $K$-theory.}
\label{sec_tekt}

Let $X$ be a compact Hausdorff space, $d \in \bN$, and $G \subseteq \ud$ a closed group. For every $\mQ \in H^1(X,QG)$, we consider the associated special \sC category $T \mQ \in $ $\tens ( X , \wa G )$, and define the abelian group
\begin{equation}
\label{def_tkt}
K_\mQ^0(X) := K_0 (T \mQ_{+,s}) \ \ ,
\end{equation}
\noindent according to Def.\ref{def_k}. $K_\mQ^0(X)$ is called the {\em twisted equivariant K-theory of} $X$. Note that we close for subobjects {\em after} having perfomed the additive completition. This implies that there is an immersion ${\bf vect}(X) \hra T \mQ_{+,s}$, in fact every $\mE \in {\bf vect}(X)$ appears as a subobject of some $\iota_n :=$ $( \iota , \ldots , \iota ) \in$ $T \mQ_+$. Thus, there is a morphism $K^0(X) \ra K_\mQ^0(X)$. 
We briefly discuss the relationship between $K_\mQ^0(X)$ and well-known $K$-theory groups.

\begin{enumerate}
\item  If $X := \left\{ x \right\}$ reduces to a single point, then the unique
       element of $\tens ( \left\{ x \right\} ,$ $ \wa G )$ is the dual $\wa G$. 
       It is well-known that by closing $\wa G$ w.r.t. direct sums and subobjects
       we get all the finite-dimensional representations of $G$. Thus, 
       $K_\mQ (\left\{ x \right\})$ coincides with (the additive group) of the
       representation ring $R(G)$.
\item  If $G = {\mathbb I}_d$ is the trivial subgroup of $\ud$ (so that
       $QG = \ud$) and $\mQ \in H^1(X,\ud)$, then $T \mQ = \wa \mE$ for some rank 
       $d$ vector bundle $\mE \ra X$. Thus, $K_\mQ^0(X)$ coincides with 
       $K^0(X)$. 
       %
       %
\item  Let $X$ be a trivial $G$-space. We consider the trivial principal 
       $QG$-bundle over $X$, say $\mQ_0$. The tensor \sC category $T \mQ_0$
       is isomorphic to the trivial special category $X \times \wa G$. Now,
       every finite-dimensional $G$-Hilbert space $M$ appears as an object
       of the category $\wa G_{+,s}$ obtained by closing $\wa G$ w.r.t. direct 
       sums and subobjects. Moreover, every $G$-vector bundle $\mE \ra X$ in 
       the sense of \cite{Seg68} is a direct summand of some trivial bundle 
       $X \times M$ (\cite[Prop.2.4]{Seg68}). Thus, $\mE \in$ ${\bf obj} \ 
       (X \times \wa G)_{+,s}$, and $K^0_{\mQ_0}(X)$ is isomorphic 
       to the equivariant $K$-theory $K_G(X)$ in the sense of \cite{Seg68}. 
       Note that for a trivial $G$-space, we find $K^0_{\mQ_0}(X) \simeq$ 
       $K_G(X) \simeq$ $R(G) \otimes K^0(X)$ (see \cite[Prop.2.2]{Seg68}.
\item  Let $\mQ \in H^1(X,QG)$ such that $\mQ = p_* \mN$ for some $\mN \in$
       $H^1(X,NG)$, $\mN :=$ 
       $( \left\{ X_i \right\} , \left\{ u_{ij} \right\} )$. 
       Then, there exists a vector bundle $\mE \ra X$ with 
       associated $\ud$-cocycle $i_{\ud,*} \mN \in$ $H^1(X,\ud)$, and a $G$-fibre 
       bundle $\mG \ra X$ with associated ${\bf aut} G$-cocycle 
       ${\bf ad}_* \mN \in$ $H^1(X,{\bf aut}G)$. By Ex.\ref{ex_d_vb}, 
       we find that $\mE$ is a $\mG$-equivariant vector bundle
       in the sense of \cite[\S 1]{NT04}. Now, $T\mQ$ is isomorphic 
       to the category $\wa \mE_\mG$
       introduced in Ex.\ref{ex_d_gb}, thus there is a morphism
       \begin{equation}
       \label{eq_mtn}
       K^0_\mQ(X) \ra K^0_\mG(X) \ \ ,
       \end{equation}
       \noindent where $K^0_\mG(X)$ denotes the gauge-equivariant $K$-theory
       of $X$ in the sense of \cite[\S 3]{NT04}. At the present moment, it is 
       not clear whether (\ref{eq_mtn}) is one-to-one. In fact, there is no 
       evidence that by 
       closing $\wa \mE_\mG$ w.r.t. direct sums and subobjects we get all the 
       $\mG$-equivariant vector bundles over $X$. Since it is possible to 
       reconstruct $\mG$ starting from the dual category $\wa \mE_\mG$ (see 
       \cite[Thm.7.3]{Vas05}), it is natural to conjecture that (\ref{eq_mtn}) 
       is an isomorphism. This point is object of a work in progress.
       %
       %
\end{enumerate}

In the case in which $\mQ$ does not belong to the image of the map (\ref{def_p*}), the group $K_\mQ^0(X)$ cannot be interpreted in terms of usual (equivariant) $K$-theory; up to direct sums, the elements of $K_\mQ^0(X)$ arise from projections belonging to the $C(X)$-algebras $\mM_{\mQ,r,r}$, $r \in \bN$ (see Sec.\ref{ssec_sc}). These $C(X)$-algebras have fibre $(H_d^r , H_d^r)_G$, anyway cannot be interpreted in terms of equivariant endomorphisms of the $r$-fold tensor power of some rank $d$ vector bundle, as claimed in Thm.\ref{thm_ci}.
Something similar happens in the setting of twisted $K$-theory considered in \cite{Ati}: given an element $\mP \in$ $H^2(X,\bT) \simeq$ $H^3(X,\bZ)$, a $K$-theory group $K_\mP^0(X)$ is constructed. The elements of $K_\mP^0(X)$ are projections of a $C(X)$-algebra $\mA_\mP$ with fibre the \sC algebra $\mK$ of compact operators; if $\mP \neq 0$, then $\mA_\mP$ is not isomorphic to $C(X) \otimes \mK$, and cannot be interpreted as the algebra of compact endomorphisms of some bundle of Hilbert spaces. $\mP$ is called the Dixmier-Douady class of $\mA_\mP$, and this is the reason why we adopted the same terminology for the invariant $\delta (\mQ)$.

Basic properties and applications of $K_\mQ^0(X)$ are objects of a work in progress. We just mention the fact the $K_\mQ^0(X)$ has a natural ring structure arising from the tensor product of $T \mQ$; such a ring structure plays an important role in the computation of the $K$-theory of the \sC algebra $\mO_\mQ$.


\end{document}